\documentclass[10pt,twoside]{siamltex}
\usepackage{amsfonts,epsfig}
\usepackage{epic,eepic}

\begin{document}

\bibliographystyle{plain}
\title{
Matrices Totally Positive Relative to
a Tree\footnote{This research was conducted during the
summer of 2007 as part of the College of William and Mary
REU program and was supported by NSF grant DMS-03-53510}}

\author{
Charles R. Johnson\thanks{Department of Mathematics,
College of William and Mary, Williamsburg, VA 23187\newline
(crjohnso@math.wm.edu).}
\and
Roberto S. Costas-Santos\thanks{Department of Mathematics,
College of William and Mary, Williamsburg, VA 23187 \newline
(rcostassantos@wm.edu).
The research of RSCS has been supported by Direcci\'{o}n General de Investigaci\'{o}n del
Ministerio de Educaci\'{o}n y Ciencia of Spain under grant MTM
2006-13000-C03-02..}
\and
Boris Tadchiev\thanks{Undergraduate, Binghamton University, Binghamton,
NY 13902 (btadchi1@binghamton.edu)}}
\pagestyle{myheadings}
\markboth{C. R. Johnson, R. S. Costas-Santos, and B. Tadchiev}
{Matrices Totally Positive Relative to a Tree}
\maketitle

\begin{abstract}
It is known that for a totally positive (TP) matrix, the eigenvalues
are positive and distinct and the eigenvector associated with the
smallest eigenvalue is totally nonzero and has an alternating sign
pattern.
Here, a certain weakening of the TP hypothesis is shown to yield a
similar conclusion.
\end{abstract}


\section{Introduction} \label{intro-sec}
A matrix is called \textit{totally positive} (TP) if \textit{every}
minor of it is positive. We will be interested in submatrices of a
given matrix that are TP, or permutation similar to TP.
Thus, we will be interested in permuted submatrices, identified by
\textit{ordered} index lists.

For any ordered index lists $\alpha$, $\beta \subseteq\{1,\dots,n\}$,
with $|\alpha|=|\beta|=k$, by $A[\alpha ; \beta]$ we mean the
$k$-by-$k$ submatrix  that lies in the rows of $A\in M_n$ indexed by
$\alpha$ and the columns indexed by $\beta$, and with the order of
the rows (resp. columns) determined by the order in $\alpha$ (resp. $\beta$).

Now, suppose that $T$ is a labelled tree on vertices $1,\ldots,n$,
and $A \in$ $M_{n}$.
If $P$ is an induced path of $T$, by $A[P]$ we mean $A[\alpha]$ in
which $\alpha$ consists of the indices of the vertices of $P$ in the
order in which they appear along $P$.
Since everything we discuss is independent of reversal of order,
there in no ambiguity.

For a given labelled tree $T$ on $n$ vertices, we say that
$A \in M_{n}(\mathbb{R})$ is T-TP if, for every path $P$ in $T$, $A[P]$ is TP.
Of course, if $T$ is a path with vertices labelled in natural order,
a T-TP matrix is, simply, TP, and if the path is labelled in some other
way, a T-TP matrix is permutation similar to a TP matrix.
If $T$ is a tree but not a path, only certain (re-ordered) proper
principal submatrices of a T-TP matrix are TP and the matrix itself
need not be permutation similar to a TP matrix.
For a T-TP matrix, properly less is required than for a TP matrix.
Also, like TP matrices, T-TP matrices are entry-wise positive.
This follows because there exists a path connecting vertices $i$ and $j$
in tree $T$, so that every entry $a_{ij}$ in the corresponding T-TP matrix
is in a submatrix that is, by definition, TP.
Since all the entries in a TP matrix are positive, $a_{ij}$ is positive for
all $i$ and $j$.
\\
\indent
There are many strong structural properties present in a TP matrix \cite{FJ}.
Among them is the fact that the eigenvalues are real, positive and distinct.
Of course, the largest one is the Perron root and its eigenvector may be
taken to be positive.
The fact that this property of a TP matrix holds for T-TP matrices is clear
from the fact that the entries are positive.
The eigenvectors of the remaining eigenvalues alternate in sign subject
to well-defined requirements, and, in particular, the eigenvector,
associated with the smallest eigenvalue, alternates in sign
as: ($+, -, +, -,\ldots,\ $).
This is 
because the inverse, or adjoint, has a checkerboard sign
pattern and the Perron root of the alternating sign signature similarity
of the inverse is the inverse of the smallest eigenvalue of the original TP matrix.
\\
\indent
If we return to the view that a TP matrix is one that is T-TP relative
to a naturally labelled path $T$, then the sign pattern of the ``last"
eigenvector may be viewed as alternation associated with each edge of
the path, i.e. if $\{i,j\}$ is an edge of $T$, then $v_{i}v_{j}<0$ for
the eigenvector $v$ associated with the smallest eigenvalue.
We may view any labelled tree as imposing a sign pattern on a totally
nonzero vector in an analogous way.
We say that the vector $v \in \mathbb{R}^{n}$ is signed according to
the labelled tree $T$ on $n$ vertices if, whenever $\{i,j\}$ is an edge
of $T$, then $v_{i}v_{j}<0$.
This means that $v$ is totally nonzero and that the sign pattern of $v$
is uniquely determined, up to a factor of $-1$.
We know that the eigenvector associated with the smallest eigenvalue of
a TP matrix is signed according to the standardly labelled path $T$
(relative to which the TP matrix is T-TP).
It is an easy exercise to see that if the path is labelled in some other
way, a T-TP matrix still has the ``last" eigenvector signed according to
the alternatively labelled path.
\\
\indent
J. Garloff \cite{G} relayed to us an old conjecture of
A. Neumaier that for any tree $T$, the eigenvector
associated with the smallest eigenvalue of a T-TP matrix
should be signed according to the labelled tree $T$ (see
note at end).
We refer to this as ``the T-TP conjecture".
We note that, again, via permutation similarity, the
labelling of the tree, per se, is not important.
If the conjecture were correct for one labelling of a given
tree, it would be correct for another.

We find here that the T-TP conjecture is correct for trees
on fewer than 5 vertices.
This requires proof only for the star on 4 vertices and this
situation is remarkably intricate, so that the proof itself
is of interest.
However, the general conjecture is false, as we see by example
for all non-paths on 5 (or more) vertices, with or without
an assumption of symmetry.
This raises the natural question of whether the basic T-TP
hypothesis may be augmented in some natural, but limited,
way to obtain the smallest eigenvector conclusion.
We find that it can be for the star on $n$ vertices.

Our arguments will heavily appeal to the adjoint of a T-TP
matrix (or one satisfying additional hypotheses), as a
surrogate for the inverse, and we frequently use Sylvester's
determinantal identity \cite{HoJo85}, along with ad hoc
arguments, to determine the sign pattern of the adjoint.
Let $\tilde{A}=Adj(A)$ for $A \in M_{n}$, so that
\begin{equation} \label{aux1}
A \tilde{A} = \tilde{A} A = (\det A)I.
\end{equation}

Note that if $\tilde{A}$ is signature similar to a positive matrix,
then $\tilde{A}$ has a Perron root and there are two possibilities:
either (1) $A$ is invertible and the smallest eigenvalue of $A$
is real and its eigenspace is that of the Perron root of $\tilde{A}$
or (2) $A$ is singular of rank $n-1$ and the nullspace of $A$ (or
eigenspace associated with $0$, the smallest eigenvalue) is spanned
by any column of $\tilde{A}$.
In either event, the smallest eigenvalue of $A$ is real, has
multiplicity one and its typical eigenvector is totally nonzero
with sign pattern corresponding to that of a column of $\tilde{A}$
or, equivalently, the diagonal of the signature matrix by which
$\tilde{A}$ is signature similar to a positive matrix.

The version of Sylvester's identity we shall often use is the following:
\begin{equation} \label{sylv-iden}
\det A[\alpha;\beta]=\frac{\det A[\alpha';\beta'] \det A['\!\alpha;'\!\!\beta] -
\det A[\alpha';'\!\!\beta] \det A['\!\alpha;\beta']}{\det A['\!\alpha';'\!\!\beta']},
\end{equation}
in which $\alpha$ and $\beta$ are index sets of the same size, $\alpha'$
(resp. $\beta'$) is $\alpha$ (resp. $\beta$) without the last index,
$'\!\alpha$ (resp. $'\!\beta$) is $\alpha$
(resp. $\beta$) without the first index, and $'\!\alpha'$ (resp.
$'\!\beta'$) is $\alpha$ (resp. $\beta$)
without the first index and last index.
We also adopt the notation that a \verb+^+ over an index in an index set
means that the index is omitted from the set. Note that, above, as
throughout, these index sets are ordered.
\section{The star on 4 vertices}
$\  $\bigskip
\begin{theorem}\label{the2.1}
For any labelled tree $T$ on fewer than 5 vertices, any T-TP matrix
has smallest eigenvalue that is real and a totally nonzero eigenvector
that is signed according to $T$.
\end{theorem}
\begin{proof}
Since the only tree on fewer than 5 vertices that is not a path
is the star on 4 vertices, we need only consider such a star.
Since the claim is independent of the particular labelling, we
need consider only the star on 4 vertices with center vertex 1.
Then we wish to show that if A is T-TP, then the sign pattern of
$\tilde{A}$ is
\begin{equation} \label{sign-pattern}
\left[
  \begin{array}{cccc}
    + & - & - & - \\
    - & + & + & + \\
    - & + & + & + \\
    - & + & + & + \\
  \end{array}
\right].
\end{equation}
Whether $\det A <0$, equal to $0$, or $>0$, this will suffice to
show both the reality of the smallest eigenvalue and that its
eigenvector is signed as
\begin{equation}\label{sign-pattern2}
\pm
\left[
  \begin{array}{c}
    +  \\
    -  \\
    -  \\
    -  \\
  \end{array}
\right].
\end{equation}
Our intent is to determine the sign of each entry in $\tilde{A}$.
Since
\begin{equation}\label{aux3}
\tilde{a}_{i,j}=(-1)^{i+j} \det A [1,\ldots,\hat{j},n; 1,\ldots,\hat{i},n],
\end{equation}
we can determine these signs by calculating each 3-by-3 minor of $A$
using only the T-TP hypothesis.
Unfortunately, this must be done in different ways for different minors.
Using Sylvester's identity (\ref{sylv-iden}), we get that
\begin{eqnarray*}
\tilde{a}_{3,2}
&=&
(-1)^{3+2} \det A [1,3,4;1,2,4]
\\
&=&
- \det A [3,1,4;2,1,4]
\\
&=&
- \frac{ \det A[3,1;2,1] \det A[1,4;1,4] - \det A[3,1;1,4]
\det A[1,4;2,1]}{\det A[1;1]}.
\end{eqnarray*}
Since $A[3,1,2;3,1,2]$ is TP by hypothesis, the submatrix $A[3,1;1,2]$ has
a positive determinant.
By permuting columns 1 and 2 in this submatrix, we get that
$\det A [3,1;2,1]<0$.
Likewise, $\det A[1,4;1,4]$, $\det A[3,1;1,4]$, $\ $and $\ \det A[1,4;2,1]$
are all positive since they are minors of submatrices that are also TP
by hypothesis.
Lastly, since T-TP matrices are entry-wise positive, $\det A[1;1]$ is
positive.
In an informal notation that should be unambiguous, we then have
\begin{equation} \label{entrie32}
\tilde{a}_{3,2}=- \frac{(-)(+)-(+)(+)}{(+)}>0.
\end{equation}
Similarly, we obtain
\begin{equation} \label{entrie42}\begin{array}{rl}
\tilde{a}_{4,2} =& \displaystyle (-1)^{4+2} \det A [1,3,4;1,2,3] =-
\det A [4,1,3;2,1,3]\\[3mm]  =& \displaystyle - \frac{\det A[4,1;2,1]
\det A[1,3;1,3] - \det A[4,1;1,3] \det A[1,3;2,1]}{\det A[1;1]}
\\ =& \displaystyle - \frac{(-)(+)-(+)(+)}{(+)}>0,
\end{array}\end{equation}
and
\begin{equation} \label{entrie43}\begin{array}{rl}
\tilde{a}_{4,3} =& \displaystyle (-1)^{4+3} \det A [1,2,4;1,2,3]
=- \det A [4,1,2;3,1,2] \\[3mm] =& \displaystyle - \frac{\det A[4,1;3,1]
\det A[1,2;1,2] - \det A[4,1;1,2] \det A[1,2;3,1]}{\det A[1;1]}
\\ =& \displaystyle - \frac{(-)(+)-(+)(+)}{(+)}>0.
\end{array} \end{equation}

By interchanging rows and columns in the above computations, we
see that $\tilde{a}_{2,3}$, $\tilde{a}_{2,4}$, and
$\tilde{a}_{3,4}$ are also positive.
To determine the signs of other entries in $\tilde{A}$, different
methods need to be used: we already proved that
$\det A[1,3,4;1,2,3]>0$ so $\det A[1,3,4;2,1,3]<0$.
Applying Sylvester's identity, we get the following inequality:
\begin{equation} \label{ineq1}
\frac{\det A[1,3;2,1] \det A[3,4;3,1]}{\det A[1,3;1,3]}
>\det A[3,4;1,2].
\end{equation}
It is possible to use this inequality to determine the
sign of $\tilde{a}_{4,1}$:
$$ \begin{array}{@{\!}r@{\!}l}
\tilde{a}_{4,1}=& \displaystyle (-1)^{4+1} \det A[2,3,4;1,2,3] =
\det A[2,3,4;2,1,3] \\[3mm] =& \displaystyle  \frac{\det A[2,3;2,1]
\det A[3,4;1,3] - \det A[2,3;1,3] \det A[3,4;2,1]}{\det A[3;1]}
\\[3mm] =& \displaystyle  \frac{- \det A[2,3;2,1] \det A[3,4;3,1] +
\det A[2,3;1,3] \det A[3,4;1,2]}{\det A[3;1]} \\[3mm]
<& \displaystyle
\frac{- \det A[2,3;2,1] \det A[3,4;3,1] + \det A[2,3;1,3]
\frac{\det A[1,3;2,1] \det A[3,4;3,1]}{\det A[1,3;1,3]}}{\det A[3;1]}
\\[4mm]
=& \displaystyle \frac{\det A[3,4;3,1]}{\det A[1,3;1,3]} \frac{\det A[1,3;2,1]
\det A[2,3;1,3] - \det A[2,3;2,1] \det A[1,3;1,3]}{\det A[3;1]}
\\[5mm]
=& \displaystyle \frac{\det A[3,4;3,1]}{\det A[1,3;1,3]} \frac{\det A[2,3;2,1]
\det A[3,1;1,3] - \det A[2,3;1,3] \det A[3,1;2,1]}{\det A[3;1]}
\\[5mm]
=& \displaystyle \frac{\det A[3,4;3,1]}{\det A[1,3;1,3]} \det A[2,3,1;2,1,3]
= - \frac{\det A[3,4;3,1]}{\det A[1,3;1,3]} \det A[2,1,3;2,1,3]<0.
\end{array}
$$
Using a similar argument gives that $\tilde{a}_{14}<0$.
Moreover, by the symmetry of the graph and by interchanging rows
and columns $2$ and $4$, it is straightforward that $A[1,4,3,2]$
is also T-TP, so $\tilde{a}_{21}$ and $\tilde{a}_{12}$ are also
negative.
Similarly, by interchanging rows and columns $3$ and $4$,
we get that $\tilde{a}_{31}$ and $\tilde{a}_{13}$ are negative
as well.
Also, it follows from the hypothesis that $\tilde{a}_{22}$,
$\tilde{a}_{33}$, and $\tilde{a}_{44}$ are positive.
Now we have that
$$
\tilde{A}=
\left[
  \begin{array}{cccc}
    \tilde{a}_{11} & - & - & - \\
    - & + & + & + \\
    - & + & + & + \\
    - & + & + & + \\
  \end{array}
\right].
$$
Taking into account (\ref{aux1}), we get that
\begin{eqnarray*}
\left[
  \begin{array}{cccc}
    \tilde{a}_{11} & - & - & - \\
    - & + & + & + \\
    - & + & + & + \\
    - & + & + & + \\
  \end{array}
\right]
\left[
  \begin{array}{cccc}
    + & + & + & + \\
    + & + & + & + \\
    + & + & + & + \\
    + & + & + & + \\
  \end{array}
\right]
&=&
\left[
  \begin{array}{cccc}
    \det A & 0 & 0 & 0 \\
    0 & \det A & 0 & 0 \\
    0 & 0 & \det A & 0 \\
    0 & 0 & 0 & \det A \\
  \end{array}
\right]
\end{eqnarray*}
and multiplying the first row of $\tilde{A}$ by the second
column of $A$ we get
$$
\tilde{a}_{11}+(-)+(-)+(-)=0,
$$
so $\tilde{a}_{11}$ has to be positive and $\tilde{A}$ has the
desired sign pattern.
\end{proof}
\section{Numerical examples}
In order to understand this phenomenon better, let us consider
an example.
The following 4-by-4 matrix is T-TP relative to the star with
center vertex 1:
\begin{eqnarray*}
A=
\left[
\begin{array}{cccc}
 130 & 78 & 98 & 96 \\
 90 & 108 & 34 & 25 \\
 116 & 57 & 137 & 44 \\
 55 & 1 & 39 & 112
\end{array}
\right].
\end{eqnarray*}
The eigenvector associated with the smallest eigenvalue, which  is
$\lambda_{4}\approx2.5$, has each entry rounded to the nearest
hundredth, in fact
\begin{eqnarray*}
\textbf{x}
\approx
\left[
\begin{array}{c}
 -3.12 \\
 1.93 \\
 1.55 \\
 1
\end{array}
\right]=
-\left[
\begin{array}{c}
 3.12 \\
 -1.93 \\
 -1.55 \\
 -1
\end{array}
\right].
\end{eqnarray*}
The adjoint of A is
\begin{eqnarray*}
\tilde{A}
=
\left[
\begin{array}{cccc}
 1308414 & -641920 & -560896 & -757860 \\
 -791797 & 446528 & 327360 & 450406 \\
 -646651 & 290240 & 328640 & 360378 \\
 -410282 & 210176 & 158080 & 292428
\end{array}
\right].
\end{eqnarray*}
Both \textbf{x} and $\tilde{A}$ have the predicted sign pattern.
Note that $A$ has 2 complex eigenvalues $\lambda= 83.6571\pm 4.24099 i$
which are greater, in absolute value, than $\lambda_4$.
In fact, it is possible to reproduce examples, for instance
with real spectrum and any sign pattern for the intermediate eigenvectors.

Next we note that non-path claims, similar to those of
theorem \ref{the2.1}, may not  be made for trees on larger
numbers of vertices.
Also, as said before, similar statements are not valid for
intermediate eigenvalues/eigenvectors.
It suffices to consider 5 vertices, and for trees on 5 vertices
there are two that are not paths: the ``star" and the ``pitchfork."

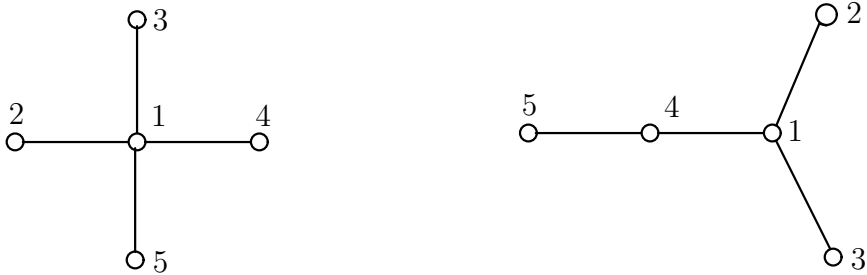
\begin{center}
\begin{minipage}{60mm}
\setlength{\unitlength}{0.00083333in}
\begingroup\makeatletter\ifx\SetFigFont\undefined%
\gdef\SetFigFont#1#2#3#4#5{%
  \reset@font\fontsize{#1}{#2pt}%
  \fontfamily{#3}\fontseries{#4}\fontshape{#5}%
  \selectfont}%
\fi\endgroup%
{\renewcommand{\dashlinestretch}{30}
\begin{picture}(1658,1727)(0,-10)
\thicklines
\put(68,812){\circle{106}}
\put(828,812){\circle{106}}
\put(1590,812){\circle{106}}
\put(816,75){\circle{106}}
\put(828,1575){\circle{106}}
\drawline(116,812)(770,812)
\drawline(881,812)(1539,812)
\drawline(828,867)(828,1535)
\drawline(816,125)(816,773)
\put(28,925){\makebox(0,0)[lb]{{\SetFigFont{12}{14.4}{\rmdefault}{\mddefault}{\updefault}2}}}
\put(928,0){\makebox(0,0)[lb]{{\SetFigFont{12}{14.4}{\rmdefault}{\mddefault}{\updefault}5}}}
\put(1566,912){\makebox(0,0)[lb]{{\SetFigFont{12}{14.4}{\rmdefault}{\mddefault}{\updefault}4}}}
\put(916,912){\makebox(0,0)[lb]{{\SetFigFont{12}{14.4}{\rmdefault}{\mddefault}{\updefault}1}}}
\put(928,1512){\makebox(0,0)[lb]{{\SetFigFont{12}{14.4}{\rmdefault}{\mddefault}{\updefault}3}}}
\end{picture}
}
\end{minipage}
\qquad
\begin{minipage}{60mm}
\setlength{\unitlength}{0.00083333in}
\begingroup\makeatletter\ifx\SetFigFont\undefined%
\gdef\SetFigFont#1#2#3#4#5{%
  \reset@font\fontsize{#1}{#2pt}%
  \fontfamily{#3}\fontseries{#4}\fontshape{#5}%
  \selectfont}%
\fi\endgroup%
{\renewcommand{\dashlinestretch}{30}
\begin{picture}(2142,1691)(0,-10)
\thicklines
\put(68,850){\circle{106}}
\put(828,850){\circle{106}}
\put(1590,850){\circle{106}}
\put(1969,75){\circle{106}}
\put(1928,1588){\circle{128}}
\drawline(116,850)(774,850)
\drawline(891,850)(1539,850)
\drawline(1616,897)(1884,1532)(1928,1521)
\drawline(1616,798)(1952,126)
\put(28,963){\makebox(0,0)[lb]{{\SetFigFont{12}{14.4}{\rmdefault}{\mddefault}{\updefault}5}}}
\put(1686,800){\makebox(0,0)[lb]{{\SetFigFont{12}{14.4}{\rmdefault}{\mddefault}{\updefault}1}}}
\put(916,950){\makebox(0,0)[lb]{{\SetFigFont{12}{14.4}{\rmdefault}{\mddefault}{\updefault}4}}}
\put(2053,1538){\makebox(0,0)[lb]{{\SetFigFont{12}{14.4}{\rmdefault}{\mddefault}{\updefault}2}}}
\put(2081,0){\makebox(0,0)[lb]{{\SetFigFont{12}{14.4}{\rmdefault}{\mddefault}{\updefault}3}}}
\end{picture}
}
\end{minipage}
\end{center}


Again the labelling is immaterial.
Let us consider the following matrix that is T-TP relative to
the star above:
\begin{eqnarray*}
A
=
\left[
  \begin{array}{ccccc}
    55 & 77 & 10 & 17 & 49 \\
    40 & 84 & 3 & 1 & 8 \\
    57 & 74 & 86 & 15 & 47 \\
    94 & 2 & 8 & 86 & 58 \\
    48 & 41 & 4 & 4 & 78 \\
  \end{array}
\right]
\end{eqnarray*}
Note that in this example, the eigenvector associated with the smallest
eigenvalue, $\lambda_{5}\approx-6.16$, does not have the predicted sign
pattern.
Here is the eigenvector in question, with each entry approximated to
the nearest hundredth:
\begin{eqnarray*}
\textbf{x}
\approx
\left[
  \begin{array}{c}
    -2.98 \\
    1.21 \\
    -0.02 \\
    2.39 \\
    1 \\
  \end{array}
\right].
\end{eqnarray*}
The sign pattern in \textbf{x} is not consistent with the conjecture
because the adjoint of A does not have the correct sign pattern
(specifically, entries $\tilde{a}_{31}$ and $\tilde{a}_{33}$ have
the ``wrong" sign), and, consequently, S$\tilde{A}$S is not positive, where
\begin{eqnarray*}
\tilde{A}
=
\left[
  \begin{array}{ccccc}
    42023084 & -27857784 & -2494736 & -6756454 & -17014640 \\
    -18274672 & 7046528 & 1241168 & 2950496 & 7815680 \\
    2070092 & 1908264 & -5017752 & 386110 & 1240248 \\
    -35907780 & 21866360 & 2481608 & 951670 & 18111768 \\
    -14519176 & 12220096 & 1012872 & 2538312 & 279496 \\
  \end{array}
\right].
\end{eqnarray*}
Let us now consider the following symmetric matrix that is
T-TP relative to the pitchfork above:
\begin{eqnarray*}
A
=
\left[
  \begin{array}{ccccc}
    88 & 50 & 35 & 78 & 38 \\
    50 & 48 & 19 & 27 & 11 \\
    35 & 19 & 41 & 13 & 6 \\
    78 & 27 & 13 & 86 & 44 \\
    38 & 11 & 6 & 44 & 59 \\
  \end{array}
\right]
\end{eqnarray*}
Here, as in the previous example, the eigenvector associated with
the smallest eigenvalue, $\lambda_{5}\approx-2.54$, does not have
the predicted sign pattern.
The following is the eigenvector in question, with each entry
approximated to the nearest hundredth:
\begin{eqnarray*}
\textbf{x}
\approx
\left[
  \begin{array}{c}
    -68.08 \\
    32.75 \\
    26.69 \\
    45.57 \\
    1 \\
  \end{array}
\right].
\end{eqnarray*}
\indent
\section{The star on $n$ vertices}
This naturally raises the question of how the T-TP hypothesis may
be augmented for larger trees to obtain the smallest
eigenvalue/eigenvector conclusion, without assuming as much as TP.
Note that for stars on larger numbers of vertices, the T-TP
hypothesis applies to only rather ``small" submatrices.
However, in this case there seems to be a natural augmentation of
the T-TP hypothesis that leads to the desired conclusion.
For the purposes of stating this additional hypothesis, we will
define ``submatrices associated with the deletion of pendant
vertices" as $A[1,\ldots,\hat{p_{k}},\ldots,n]$, where $p_{k}$
is a pendant (degree 1) vertex of tree $T$ on $n$ vertices.
The additional hypothesis is that these submatrices must be
$P$-matrices, or matrices that have only positive principal minors.
Note that in the case of the star on 4 vertices this adds
nothing to the T-TP hypothesis.

\begin{theorem}  \label{theo41}
Let $T$ be a star on $n$ vertices. Suppose that $A$ is T-TP and
that all the submatrices of $A$ associated with the deletion of
pendant vertices are $P$-matrices.
Then, the smallest eigenvalue of $A$ is real, has multiplicity
one and has an eigenvector signed according to $T$.
\end{theorem}
\begin{proof}
Let us consider the following star $T$:

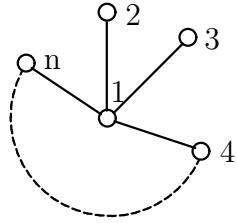
\begin{center}
\setlength{\unitlength}{0.00083333in}
\begingroup\makeatletter\ifx\SetFigFont\undefined%
\gdef\SetFigFont#1#2#3#4#5{%
  \reset@font\fontsize{#1}{#2pt}%
  \fontfamily{#3}\fontseries{#4}\fontshape{#5}%
  \selectfont}%
\fi\endgroup%
{\renewcommand{\dashlinestretch}{30}
\begin{picture}(1415,1489)(0,-10)
\thicklines
\drawline(616,728)(616,1280)
\put(620,682){\circle{106}}
\put(1204,476){\circle{106}}
\put(112,1025){\circle{106}}
\put(618,1343){\circle{106}}
\put(1122,1187){\circle{106}}
\drawline(148,994)(580,706)
\put(1226,1116){\makebox(0,0)[lb]{{\SetFigFont{12}{14.4}{\rmdefault}{\mddefault}{\updefault}3}}}
\drawline(664,662)(1168,494)
\thinlines
\thicklines
\drawline(659,714)(1091,1146)
\put(640,781){\makebox(0,0)[lb]{{\SetFigFont{12}{14.4}{\rmdefault}{\mddefault}{\updefault}1}}}
\put(1321,407){\makebox(0,0)[lb]{{\SetFigFont{12}{14.4}{\rmdefault}{\mddefault}{\updefault}4}}}
\put(222,988){\makebox(0,0)[lb]{{\SetFigFont{12}{14.4}{\rmdefault}{\mddefault}{\updefault}n}}}
\put(732,1267){\makebox(0,0)[lb]{{\SetFigFont{12}{14.4}{\rmdefault}{\mddefault}{\updefault}2}}}
\dashline{60.000}(87.000,976.000)(56.073,907.253)(33.748,835.252)
	(20.357,761.068)(16.098,685.805)(21.036,610.584)
	(35.097,536.524)(58.072,464.727)(89.618,396.262)
	(129.266,332.148)(176.427,273.339)(230.397,220.710)
	(290.375,175.045)(355.467,137.023)(424.704,107.211)
	(497.057,86.051)(571.447,73.859)(646.769,70.816)
	(721.900,76.968)(795.724,92.222)(867.140,116.353)
	(935.087,149.000)(998.553,189.679)(1056.593,237.783)
	(1108.343,292.596)(1153.034,353.303)(1190.000,419.000)
\end{picture}
}
\end{center}

Again, as in the 4-by-4 case, if the theorem can be proven for
a particular labelling of the tree, then it will be true for
all other labellings.
In order for an $n$-by-$n$ T-TP matrix to have the correct
sign pattern in the ``last" eigenvector, $\tilde{A}$ has to
have the following sign pattern:
\begin{eqnarray*}
\tilde{A}
=
\left[
  \begin{array}{c|ccc}
    + & - & \cdots & - \\\hline
    - & + & \cdots & + \\
    \vdots & \vdots & \ddots & \vdots \\
    - & + & \cdots & + \\
  \end{array}
\right].
\end{eqnarray*}
Let us look at permuted minors of the form
$\det A[i,1,\aleph;j,1,\aleph]$, where $\aleph$ is an ordered index
list, in lexicographic order, not containing $i$, $j$, and 1, such
that $i \ne j$ and $i \ne 1 \ne j$. Also, let $\aleph'$ be $\aleph$
without the last entry $l$.
In the 4-by-4 case, these are the minors $\det A[2,1,4;3,1,4]$,
$\det A[2,1,3;4,1,3]$, and $\det A[3,1,2;4,1,2]$, as well as the
minors obtained by replacing rows by columns.
It was already shown what signs the permutations of these minors
must have, and it follows that the above are negative.
Using induction and assuming that the minors of this form are
also negative in the $n$-by-$n$ case, we have the following in
the $(n+1)$-by-$(n+1)$ case:
$$\begin{array}{rl}
\det &\hspace{-3mm}A[i,1,\aleph;j,1,\aleph]=\\[3mm]
& \displaystyle \frac{\det A[i,1,\aleph';j,1,\aleph']
\det A[1,\aleph;1,\aleph] - \det A[i,1,\aleph';1,\aleph]
\det A[1,\aleph;j,1,\aleph']}{\det A[1,\aleph';1,\aleph']}.
\end{array}
$$
In the third minor in the numerator, we can move the last element
$l$ of $\aleph$ to the beginning of the list, resulting in
$\det A[i,1,\aleph';l,1,\aleph']$.
The same can be done to the fourth minor in the numerator, resulting
in $\det A[l,1,\aleph';j,1,\aleph']$. Since the number of changes
in both instances is the same, the total number of these permutations
is even, and there is no minus sign introduced.
Thus, the above can be rewritten as follows:
$$ \begin{array}{rl}
\det &\hspace{-3mm} A[i,1,\aleph;j,1,\aleph]=\\[3mm]
& \displaystyle \hspace{-2mm}\frac{\det A[i,1,\aleph';j,1,\aleph']
\det A[1,\aleph;1,\aleph]\! -\! \det A[i,1,\aleph';l,1,\aleph'] \det
A[l,1,\aleph';j,1,\aleph']}{\det A[1,\aleph';1,\aleph']}.
\end{array}$$
The first, third, and fourth minors in the numerator are negative by induction.
Now, making use of our additional hypothesis, we can determine the sign
of the remaining two minors.
This condition that the submatrices associated with the deletion
of pendant vertices are $P$-matrices simply means, in the case in which
the center vertex of the star is 1, that all principal minors up
through those of size  $(n-1)$-by-$(n-1)$ which include the first
row and column, are positive.
But these minors are exactly those remaining two minors in the above
expression: the second minor in the numerator and the minor in
the denominator.
We now get the following:
$$
\det A[i,1,\aleph;j,1,\aleph]=\frac{(-)(+)-(-)(-)}{(+)} <0.
$$
Let $\aleph_{i}$ be $\aleph$ with $i$ and let $\aleph_{j}$ be
$\aleph$ with $j$, both in numerical order.
We know that $i\neq j$, so either $j<i$ or $i<j$.
Consider the case with $j<i$:
$$
\tilde{a}_{ij} = (-1)^{i+j} \det A[1,\aleph_{i};1,\aleph_{j}].
$$
Moving $j$ to the beginning of the list results in $(-1)^{j-1}$ minus
signs being introduced.
However, since $j<i$ and $\aleph_{i}$ is in numerical order and
does not include $j$, moving $i$ to the beginning of the list
results in only $(-1)^{i-2}$ minus signs.
We now have the following:
$$\begin{array}{rl}
\tilde{a}_{ij} &=\displaystyle (-1)^{i+j}
(-1)^{j-1} (-1)^{i-2} \det A[i,1,\aleph;j,1,\aleph]\\
 &= (-1)^{2(i+j-1)-1} \det A[i,1,\aleph;j,1,\aleph]>0.
\end{array}
$$
The case where $i<j$ can be proven similarly.
This proves the desired sign pattern in all the off-diagonal
entries in the lower right $(n-1)$-by-$(n-1)$ submatrix of $\tilde{A}$.
Furthermore, it follows from our additional hypothesis that
all the $(n-1)$-by-$(n-1)$ principal minors of $A$, with the
exception of the bottom right one, are positive, and so all the
diagonal entries of $\tilde{A}$, with the exception of
$\tilde{a}_{11}$, are positive.
We now have the following sign pattern for $\tilde{A}$:
\begin{eqnarray*}
\tilde{A}
=
\left[
  \begin{array}{c|ccc}
    \tilde{a}_{11} & \tilde{a}_{12} & \cdots & \tilde{a}_{1n} \\\hline
    \tilde{a}_{21} & + & \cdots & + \\
    \vdots & \vdots & \ddots & \vdots \\
    \tilde{a}_{n1} & + & \cdots & + \\
  \end{array}
\right].
\end{eqnarray*}
Now, as was previously done in the 4-by-4 case, and since
$A\tilde{A}=\tilde{A}A=\det A I$, multiplying rows $i$ of
$\tilde{A}$ by columns $j$ of $A$ and rows $k$ of $A$ by
columns $l$  of $\tilde{A}$, where $i \ne j$, $i \ne 1
\ne j$, $k \ne l$, and $k \ne 1 \ne l$, has to yield zero.
Thus, the off-diagonal entries in row and column 1 have to
be negative.
Using the same method and multiplying row 1 of $\tilde{A}$
by column 2 of $A$, we get that $\tilde{a}_{11}$ is also
positive, since that product must also result in zero.
This completes the proof.
\end{proof}
\section{Additional comments}
For trees on $n$ vertices that are not stars, examples show that
the hypothesis of theorem \ref{theo41} is not sufficient to support the
conclusion about the signing of the ``smallest" eigenvector.
It is an interesting question if there is a natural hypothesis,
more limited than TP that will.
In both theorems 2.1 and 4.1 the T-TP assumption may be broadened
slightly to ``T-oscillatory".



%
%
%
%

\bigskip
\noindent
{\bf Authors' Note:}\\
According to Garloff, Neumaier was motivated by a manuscript of Godsil,
who defined the notion of sign change of vector components, relative to a
tree, used herein.
Neumaier originally conjetured that all eigenvectors should be signed as those
of a TP matrix.
However, even the the addition of an assumption of symmetry to insure that
all eigenvalues are real, as first suggested by Garloff, does not imply that
the eigenvectors of the intermediate eigenvalues are properly signed.
This is easily seen via experimentation with examples as noted by Garloff, and
us; also the intermediate eigenvalues need not to be real without symmetry.
Since a T-TP matrix is positive, the spectral radius is an eigenvalue
and its eigenvector is uniformly signed, as in the TP case, but this
is a simple observation.
This leaves the question of the smallest eigenvalue and its eigenvector,
as raised by Neumaier and Garloff and studied here.
For this question, an assumption of symmetry  does not matter.
Though the Neumaier/Garloff conjecture (which was around for sometime)
is not generally correct, we suspect that there is an assumption, weaker
than TP, that will give  the conjectured conclusion about the smallest
eigenvalue and eigenvector.

\end{document}